\documentclass[letterpaper, 10 pt, conference]{ieeeconf}  
        
\IEEEoverridecommandlockouts                              
\usepackage{amsmath}
\usepackage{amssymb}
\usepackage[utf8]{inputenc}
\usepackage{eurosym}
\usepackage{graphicx}
\usepackage{hyperref}
\usepackage{dsfont}
\usepackage{dsfont}
\usepackage{xcolor,colortbl}
\usepackage{comment}
\usepackage{tikz}
\usepackage{subcaption}

\allowdisplaybreaks

\usepackage{hyperref}

\usepackage{ifthen}

\makeindex

\newcommand{\Rr}{{\mathbb{R}}}

\newcommand{\Ff}{{\mathcal{F}}}

\newcommand{\bx}{{\bf x}}

\newcommand{\bdx}{\dot{\bf x}}

\newcommand{\bv}{{\bf v}}

\def\dx{{\rm d}x}
\def\dy{{\rm d}y}

\def\dt{{\rm d}t}
\def\ds{{\rm d}s}
\def\leq{\leqslant}
\def\geq{\geqslant}

\newtheorem{theorem}{Theorem}
\newtheorem{problem}{Problem}
\newtheorem{rem}{Remark}

\title{\LARGE \bf
The Potential Method For Price-Formation Models
}

\author{Yuri Ashrafyan, Tigran Bakaryan, Diogo Gomes, and  Julian Gutierrez%
	\thanks{King Abdullah University of Science and Technology (KAUST), CEMSE Division, Thuwal 23955-6900. Saudi Arabia.}
		\thanks{	Email:   yuri.ashrafyan@kaust.edu.sa}
	\thanks{Email:  tigran.bakaryan@kaust.edu.sa}
	\thanks{Email:  diogo.gomes@kaust.edu.sa}
	\thanks{Email:  julian.gutierrezpineda@kaust.edu.sa}
}

\begin{document}

\maketitle
\thispagestyle{empty}
\pagestyle{empty}

\begin{abstract}
	
	We consider the mean-field game price formation model introduced by Gomes and Sa\'ude. 
	In this MFG model, agents trade a commodity whose supply can be deterministic or stochastic. Agents maximize profit, taking into account current and future prices. 
	The balance between supply and demand determines the price. We introduce
	a potential function that converts the MFG into a convex variational problem. This variational formulation 
	is particularly suitable for machine learning approaches. Here, we use a recurrent neural network
	to solve this problem.  In the last section of the paper, we compare our results
	with known analytical solutions. 
	
\end{abstract}

%

\section{Introduction}
Mean-field games (MFG) were introduced  in \cite{ll1}, \cite{Caines1} as models for systems with a large number of competitive rational agents.  In these games,  agents take decisions according to preferences  given by a utility function. Individually, agents have a negligible effect on the entire system but are affected by aggregate effects. These games arise to model power networks, free markets, social media, or pedestrians flows.

Here, we study the  electricity market price formation model from \cite{gomes2018mean}. Note, however, that the same model
describes price formation for any commodity that can be traded, stored, and whose
price is determined by market clearing conditions. 
In this model, agents have storage devices and may
charge them (buy) or supply (sell) the grid with electricity. Agents 
maximize profit by trading electricity at a price $\varpi(t)$. The balance between 
supply, $Q(t)$, and aggregate demand
determines the price, $\varpi(t)$.
We consider both deterministic and stochastic supply. 

\subsection{Deterministic supply}

The deterministic supply,  discussed in Section \ref{Sec:3-rel-det}, is 
outlined in the next problem. 
\begin{problem} \label{PMFG} 
	Let $\Theta=[0,T]\times\Rr$. 
	Suppose that $m_0\in \mathcal{P}(\Rr)$, $H\in C^2(\Rr^2)$ is uniformly convex, and $Q$ and $u_T$ are Lipschitz continuous.  
	Find $u,m: \Theta \to \Rr$ and $\varpi:[0,T]\to \Rr$ satisfying
	\begin{equation}\label{eq:MFG system}
		\begin{cases}
			-u_t +H(x,\varpi + u_x )=0 
			\\
			m_t - \left(H_p(x,\varpi + u_x)m\right)_x =0 
			\\
			-\int_{\Rr} H_p(x,\varpi + u_x)m \dx = Q(t),
		\end{cases}
	\end{equation}
	$m\geq 0$, and the initial-terminal conditions
	\begin{equation}\label{boundaryP}
		m(0,x)=m_0(x),   \quad u(T,x)=u_T(x).
	\end{equation} 
\end{problem}  

\smallskip

Existence and uniqueness of solutions was addressed in \cite{gomes2018mean}.  \cite{YTDJduality2021} examined the connection between Aubry-Mather theory and Problem \ref{PMFG}.

In 
\cite{YTDJpotential2022}, following the ideas in  \cite{DRT2021Potential}, we introduce a potential function $\varphi:\Theta\to \Rr$ giving for each $t$ the cumulative distribution of $m$.
Further, $\varphi_x=m$ and $-\varphi_t$ is the agents' current or flow. There, we formulate a variational problem and show that  \eqref{eq:MFG system} is the corresponding Euler-Lagrange equation, see Section \ref{Sec:3-rel-det}. To detail this problem, we  
let $L$ be the Legendre transform of  $H$, 
$L(x,v)=\sup_{p \in \Rr}  \left[ - p v - H(x,p) \right]$. We assume that $L$ and $H$ are uniformly convex in the second variable and that, for some $c, p>0$,
$L(x,v) \geq c|v|^p$.
Further, let $F:\Rr^2\times\Rr_0^+\to \Rr_0^+$ be the perspective function of  $L$
\begin{equation*}
	F(x,j,m) = \begin{cases} L\left(x, \frac{j}{m}\right) m, &m>0
		\\
		+\infty, & j\neq 0,\, m=0
		\\
		0, & j=0, \,m= 0. \end{cases}
\end{equation*}
The constrained variational problem corresponding to Problem \ref{PMFG} is as follows. 
\begin{problem}\label{PMFG-Var} Consider the setting of Problem \ref{PMFG}. Find  $\varphi:\Theta\to \Rr$ minimizing 
	\begin{equation*}
		\int_{\Theta} F(x,-\varphi_t,\varphi_x)-u^\prime_T(x) \varphi_t~\dx\dt, 
	\end{equation*}
	over all functions s.t. $\varphi(0,x)=\int_{-\infty}^x m_0(y)\dy$,
	and, for all $t\in [0,T]$,
	$\varphi_x(t,\cdot) \in \mathcal{P}(\Rr)$ and
	\[
	\int_{\Rr} \varphi(t,x) - \varphi(0,x) \dx = -\int_0^t Q(s)\ds.
	\]
\end{problem}

\smallskip

%
%
The functional to be minimized in the prior problem represents an average social cost; the constraint reflects the balance condition. 
In \cite{YTDJpotential2022}, we analyzed Problem \ref{PMFG-Var}. There, 
we obtained the correspondence between Problems \ref{PMFG} and \ref{PMFG-Var} and proved the existence of the price, $\varpi$, as a Lagrange multiplier.



\subsection{Stochastic supply}

The stochastic supply case was studied in \cite{GoGuRi2021} for the continuum, and in \cite{GoGuRi2021B} for $N$ players;  the authors obtained a price as a Lagrange multiplier for the balance constraint. 
In the stochastic case, discussed in Section \ref{Sec:2-mot-Stoch}, the corresponding
problem is as follows.  
\begin{problem} \label{PMFG-Stoch} Consider the setting of Problem \ref{PMFG}. Let $(\Omega, \Ff, P)$ be  a probability space that supports a standard Brownian motion $W_t$. Let $\Ff_t$ be the filtration generated by $W_t$ and let $Q$ be a progressively measurable process (PMP) with respect to $\Ff_t$. Find $m: \Theta \to \Rr$ and $u$, $Z$ and $\varpi$,  PMP with respect to $\Ff_t$, satisfying  $m\geq 0$, \eqref{boundaryP},  and
	\begin{equation}\label{eq:MFG system-Stoch}
		\begin{cases}
			-du +H(x,\varpi + u_x )dt=Z(x,t)dW_t 
			\\
			m_t - \left(H_p(x,\varpi + u_x)m\right)_x =0 
			\\
			-\int_{\Rr} H_p(x,\varpi + u_x)m \dx = Q(t). 
		\end{cases}
	\end{equation}
\end{problem} 

\smallskip

Next, in Section \ref{Sec:3-rel-Stoch},  we introduce a progressively measurable potential, $\varphi$ and derive the following stochastic variational problem associated with Problem \ref{PMFG-Stoch}.

\begin{problem}\label{PMFG-Var-Stoch} Consider the setting of Problem \ref{PMFG-Stoch}. Find a  progressively measurable  function $\varphi:\Theta\times \Omega \to \Rr$ minimizing 
	\begin{equation*}
		E\int_{\Theta} F(x,-\varphi_t,\varphi_x)   -\varpi\left(  \varphi_t+Q\varphi_x \right)-u^\prime_T \varphi_t~\dx\dt   
	\end{equation*}
	over the set of PMP  $\varphi$ 
	such that $\varphi(0,x)=\int_{-\infty}^x m_0(y)\dy$,
	and  for all $t\in [0,T]$ satisfy
	$\varphi_x(t,\cdot) \in \mathcal{P}(\Rr) $ and
	\[
	\int_{\Rr} \varphi(t,x) - \varphi(0,x) \dx = -\int_0^t Q(s)\ds.
	\]	
\end{problem}

\smallskip

In section \ref{mlern}, we
parametrize the potential by a recurrent neural network (RNN) and use the variational problem as the residual. RNN preserve progressive measurability and give a compact way to 
represent the potential. This  is particularly relevant in the random case where the problem is infinite-dimensional.  
In Section \ref{numres}, we illustrate our approach for the linear-quadratic setting and take as benchmarks the explicit formulas provided in \cite{gomes2018mean} and \cite{GoGuRi2021}.

\subsection{Prior work}

In \cite{ATM19} and \cite{alasseur2021mfg}, MFG models were used to describe intraday electricity markets. There, the price is a function of the demand rather than being determined by market clearance.  In  \cite{FTT20}, the authors considered a Stackelberg game of an intraday market, where a major agent faces many small agents.  A similar problem with  a major player was studied in  \cite{fujii2021equilibrium}. Using stochastic control theory,  in \cite{feron2021price}, the authors studied renewable energy markets in the $N$-agent and MFG settings. 
A deterministic  $N$-agent price model was considered in \cite{SummerCamp2019} and in \cite{GoGuRi2021B}. Also, 
a $N$-agent price  model  was examined in \cite{aid2020equilibrium}. There, the price is determined by the equilibrium of the system; the agents choose optimal controls (production and trading rates) to meet a demand with noise. In \cite{BS02}, \cite{BS10}, Stackelberg games modeled price formation under revenue optimization.  To study price dynamics in electricity market  in  \cite{TBD20}, the  authors considered  Cournot model and showed that their model is a MFG problem with common noise and a jump-diffusion.   The paper \cite{AidDumitrescuTankov2021} modeled the transition to renewable energies by an optimal switching MFG. 
The works \cite{JSF20} and \cite{FT20} used market-clearing conditions  to study Solar Energy Certificate Markets and flows in exchange markets, respectively.



Because there are only a few MFG models with explicit  solutions, numerical analysis of MFGs plays a crucial role, see the survey \cite{MR4214773} and the earlier works
in \cite{CDY} and \cite{DY}. 
However, prior methods cannot be used directly for our model because
the price, $\varpi$, must satisfy an integral constraint.  In \cite{YTDJpotential2022}, we introduced a numerical scheme for Problem \ref{PMFG}, which solves a convex minimization problem. In the case
of common noise, this matter is more delicate because the state space becomes infinite-dimensional. A state-space reduction strategy was developed in \cite{GoGuRi2021}
and the $N$ player case was studied in \cite{GoGuRi2021B}, but none of the prior references addressed this problem in full generality. 

\subsection{Main contributions}

This paper contains the following contributions. We present a novel formulation 
of the price problem with a random supply as a stochastic MFG system (Problem \ref{PMFG-Stoch}), extending prior works limited to the linear-quadratic and $N$-player cases. We developed
a potential approach for the  random supply problem by introducing a stochastic variational problem (Problem \ref{PMFG-Var-Stoch}). Finally, we used RNN to solve both the deterministic and stochastic price problems (Problems
\ref{PMFG-Var} and  \ref{PMFG-Var-Stoch}). Because the stochastic case is infinite-dimensional, 
the use of RNN seems to beat the curse of dimensionality. 
%
%
%


\section{Control problem associated with the price model}

Now, we describe the price formation model.

\subsection{Deterministic price problem derivation}

Consider a deterministic supply $Q$. Given the price $\varpi$,   
each agent chooses a control $v:[t,T]\to \Rr$ that minimizes
\[
u(x, t)=\inf_v\int_t^T L(\bx, v) +\varpi v +u_T(\bx(T)).
\]
By the HJ verification theorem, if $u$ solves the first equation in \eqref{eq:MFG system}, 
$u$ is the value function and $v^*=-H_p(x,\varpi + u_x )$ an optimal control. 
Because agents are rational, they use this optimal strategy; hence, 
$\bdx=v^*(t)$. Thus, the density $m$ solves the transport equation, 
the second equation in \eqref{eq:MFG system}.
While we assumed the price to be known, this MFG is a fixed point: the price must 
be consistent with 
the last identity in \eqref{eq:MFG system}
that requires aggregated total demand (left-hand side)
to match supply (right-hand side).

\subsection{Stochastic price problem derivation}
\label{Sec:2-mot-Stoch}

Consider the setting of Problem \ref{PMFG-Stoch}. 
Assume that
$Q$ and $\varpi$ are PMP with respect to $\Ff_t$. Let
$E_t[Z]=E[Z|\Ff_t]$, for any random variable $Z$.
Each agent chooses a PMP $\bv$ to minimize
\begin{equation}
	\label{vf}
	u(x, t)=E_t\int_t^T L(\bx, \bv) +\varpi \bv +u_T(\bx(T)),
\end{equation}
where $\bdx= \bv$ and $\bx(t)=x$. 
\begin{theorem}
	Let $u$ and $Z$ be PMP solving the 
	stochastic HJ equation, the first equation in \eqref{eq:MFG system-Stoch}. 
	Then, $u$ is the value function and $v=-H_p(x,\varpi + u_x )$ is an optimal control. 
\end{theorem}
\begin{rem}
	Note that the unknowns in the HJ equation are both $u$ and $Z$; the additional unknown $Z$ ensures progressive measurability. 
\end{rem}
\begin{proof}
	We have
	\begin{align*}
		&E_t\int_t^T L(\bx, \bv) +\varpi \bv ds +u_T(\bx(T))\\
		&=
		u(x,t)+E_t\int_t^T L(\bx, \bv) +\varpi v ds +du \\
		&\geq u(x,t),
	\end{align*}
	by the first equation in \eqref{eq:MFG system-Stoch} and
	the definition of $L$. Because equality holds when 
	$v=-H_p(x,\varpi + u_x )$, the claim follows. 
\end{proof}

Because agents are rational, 
each  agent evolves according to $\bdx=-H_p(\bx,\varpi + u_x )$.
This is an ODE with a random vector field. 
Accordingy, $m$ solves the (random) transport equation, 
the second equation in \eqref{eq:MFG system-Stoch}. 
Finally,  we require the balance condition in \eqref{eq:MFG system-Stoch} to hold.


\section{Potential Approach} 
\label{sec:potential_app}

Here, we describe the potential approach for our models. This approach  was introduced in \cite{DRT2021Potential}, where  one-dimensional first- and second-order MFG planning problems were examined.

\subsection{Deterministic price} 
\label{Sec:3-rel-det}

We begin our analysis by addressing the transport equation in \eqref{eq:MFG system}. 
Let $(u, m, \varpi)$ solve \eqref{eq:MFG system} with $m>0$. 
The transport equation can be
written as 
$
\text{div}_{(t,x)}\left( m , - H_p(x,\varpi + u_x)m \right) = 0.
$
Hence, by Poincar{\'e} lemma (see \cite{Csato2011ThePE}, Theorem 1.22),  there exists $\varphi : \Theta \to \Rr$ such that 
\begin{equation}\label{eq: potential relations}
	m=	\varphi_x, \quad H_p(x,\varpi + u_x)m=\varphi_t.
\end{equation}
Differentiating the Hamilton-Jacobi equation in  \eqref{eq:MFG system} with respect to $x$ and using \eqref{eq: potential relations}, we rewrite
the problem  \eqref{eq:MFG system}-\eqref{boundaryP} in terms of $\varphi$ 
\begin{equation}\label{eq: Euler-Lagrange wrt potential}
	\begin{cases}
		\left(H\left(x,-D_vL \left(x, -\tfrac{\varphi_t }{\varphi_x}\right)\right)\right)_x \\
		\qquad +\left(D_v L \left(x, -\tfrac{\varphi_t }{\varphi_x}\right)+\varpi \right)_t=0 ,
		\\
		-\int_{\Rr}  \varphi_t + Q \varphi_x ~ \dx =0,
	\end{cases}
\end{equation}
with $\varphi_x(0,\cdot)=m_0(\cdot)$ and 
\begin{equation}\label{boundary-in-phi}
	-D_vL\left(x, -\frac{\varphi_t(T,x) }{\varphi_x(T,x)}\right)-\varpi(T)=u^\prime_T(x).
\end{equation}
Simple computations show that \eqref{eq: Euler-Lagrange wrt potential} and \eqref{boundary-in-phi} are the  Euler-Lagrange equations for the functional
\begin{equation*}
	\begin{split}
		\int_{\Theta} F(x,-\varphi_t,\varphi_x)   -\varpi\left(  \varphi_t+Q\varphi_x \right)-u^\prime_T \varphi_t ~\dx\dt.   
	\end{split}
\end{equation*}
$ \varpi$ can be seen as a Lagrange multiplier for the constraint $\int_{\Rr}  \varphi_t + Q \varphi_x ~ \dx=0$. Thus, we obtain Problem \ref{PMFG-Var}. 
%
%
%
Further, 
let $\varphi\in C^2(\Theta)$ solves Problem \ref{PMFG-Var}. Then, we recover  the solution $(u, m, \varpi)$ to Problem \ref{PMFG} as follows. We set 
$
m=\varphi_x, 
$
and $u(t,x)$ is given by 
\[	
u_T(x) - \int_t^T H\left(x,-D_vL\left(x,\tfrac{-\varphi_t(s,x)}{\varphi_x(s,x)}\right)\right)\ds.
\]
$\varpi$ is given by the first equation of \eqref{eq: Euler-Lagrange wrt potential} and \eqref{boundary-in-phi}.

\subsection{Fundamental lemmas of the calculus of variations in random environments}

Before we discuss the random supply case, we need some preliminary definitions. We consider the following spaces of functions $L^2=L^2(\Theta\times\Omega)$, which consists of all PMP which are square-integrable in $\Theta\times\Omega$ and $H^1=H^1(\Theta\times\Omega)$, 
which contains all PMP, $\phi\in L^2$, such that $\phi_t, \phi_x\in L^2$. 
Note for a processe $\psi\in H^1$, $d\psi=\psi_t dt$; that is, its quadratic variation vanishes. 
We also consider the space $\tilde H(\Theta\times\Omega)$ of processes
$\psi\in L^2$ such that its stochastic differential $d\psi= a dt +b dW_t$ exists with $a, b\in L^2$.
Let $\psi\in \tilde H$ and $\phi\in H^1$. Then, we have the following identity
\begin{align}
	\label{intbp}
	&E\int_\Theta \psi \phi_t=-E\int_\Theta d(\psi) \phi
\end{align}
if $\phi(x,T)=\phi(x,0)=0$.
Moreover, if  $a, b\in L^2$ and we have 
\begin{eqnarray}
	E\int_\Theta  \varphi (a dt +b dW_t)=0, 
\end{eqnarray}
for all $\varphi\in L^2$ then $a=0$. 

\subsection{Common noise case}
\label{Sec:3-rel-Stoch}

Sufficiently regular critical points of functional in Problem \ref{PMFG-Var-Stoch} satisfy the weak form of
the Euler-Lagrange equation
\[
E\int_{\Theta} -F_j \phi_t+F_m \phi_x-\varpi \phi_t dx dt=0,
\]
for all $\phi\in H^1$.  Assume that $F_j, \varpi\in \tilde H$.
By combining the integration by parts formula and the fundamental lemma, we get
\begin{equation}\label{EL-Stoch}
	-d (-F_j-\varpi)-(F_m)_x dt =R dW_t.
\end{equation}
Next, we define $m$ and $u$ according to \eqref{eq: potential relations}, and set $R=Z_x$.  
Because of  \eqref{eq: potential relations}, we have
\begin{equation*}
	\varpi + u_x=-D_vL\big(x,-\tfrac{\varphi_t}{\varphi_x}\big)=-F_j\big(x,-\varphi_t,\varphi_x\big).
\end{equation*}
Further,  the definition of $F$, yields
\begin{align*}
	&F_m(x, -\varphi_t, \varphi_x)\\&=L\big(x,-\tfrac{\varphi_t}{\varphi_x}\big)+D_vL\big(x,-\tfrac{\varphi_t}{\varphi_x}\big)
	\tfrac{\varphi_t}{\varphi_x}\\&=
	-H(x, 	 \varpi + u_x).
\end{align*}
Therefore \eqref{EL-Stoch} can be rewritten as
\begin{equation*}
	\begin{split}
		-d(u_x) +(H(x,\varpi + u_x ))_xdt
		=RdW_t, 
	\end{split}
\end{equation*}
which is the (weak) derivative of the stochastic HJ equation in \eqref{eq:MFG system-Stoch}.

%
%
%
%


\section{Machine learning architecture}
\label{mlern}
\begin{figure}[!htb]
	\centering
	\scalebox{0.75}{
		\begin{tikzpicture}[
			roundnode/.style={circle, draw=gray!60, fill=gray!5, very thick, minimum size=6mm},
			squarednode/.style={rectangle, draw=gray!60, fill=gray!5, very thick, dashed},
			]
			\node[squarednode] at (-0.2,-1.5) [draw=none,fill=none](0Node) {} ;	
			\node[squarednode] at (1,-1.5) [draw=none,fill=none](OutputL1at0) {${\color{blue} \mathrm{h}^{\langle k-1 \rangle}}$} ;
			\draw[->,draw=blue,line width=0.25mm] (0Node.east) |- ++(0,0) |- (OutputL1at0.west);	
			
			\node[squarednode] at (4.2, -1.5) [fill=none,minimum width=50mm,minimum height=12mm,solid](RNNCellL1) {};	\node[squarednode] at (4.2, 0.2) [fill=none,minimum width=50mm,minimum height=20mm,solid](RNNCelL2) {};
			\node[squarednode] at (4.2, -0.45) [fill=none,minimum width=52mm,minimum height=36mm,solid](RNNCelL3) {};		
			\node[squarednode] at (4.2,-3.2) [draw=none,fill=none](Input1) {$\begin{array}{c} 
					\left(x^{(i)}, \mathrm{y}^{\langle k \rangle}\right)
					\\
					\mathrm{y}^{\langle k \rangle} = \left(t^{\langle k \rangle},Q^{\langle k \rangle},\mathrm{h}^{\langle k-1 \rangle}\right)
				\end{array}$} ;
			\draw[->,draw=blue,line width=0.25mm] (OutputL1at0.south) |- ++(0,0) |- (Input1.west);

			\node[squarednode] at (4.2,-1.5) [draw=none,fill=none](OutputL1at1) {$\begin{array}{c}
					{\color{blue}\mathrm{h}^{\langle k \rangle} = \tanh \left( \mathrm{W}_h \mathrm{y}^{\langle k \rangle}+ \mathrm{b}_h\right)}
				\end{array}$} ;
			\draw[->,draw=gray] (Input1.north).. controls +(up:0mm) and +(down:0mm) .. (OutputL1at1.south);				
			
			\node[squarednode] at (4.2,0.3) [draw=none,fill=none](OutputL2at1) {$\begin{array}{c}
					\mathrm{z}^{[3]}= \mathrm{W}^{[3]} \mathrm{z}^{[2]}+ \mathrm{b}^{[3]}
					\\
					\mathrm{z}^{[2]}= \mathrm{S}(\mathrm{W}^{[2]} \mathrm{z}^{[1]}+ \mathrm{b}^{[2]})
					\\
					\mathrm{z}^{[1]}=\mathrm{S} ( \mathrm{W}^{[1]} \left( x^{(i)},\mathrm{h}^{\langle k \rangle }\right) + \mathrm{b}^{[1]})
				\end{array}$} ;
			\draw[->,draw=gray] (OutputL1at1.north).. controls +(up:0mm) and +(down:0mm) .. (OutputL2at1.south);

			\node[squarednode] at (4.2,2.1) [draw=none,fill=none](OutputVat1) {$\varphi (t^{\langle k \rangle},x^{(i)}) =  \mathrm{z}^{[3]}$} ;
			\draw[->,draw=gray] (OutputL2at1.north).. controls +(up:0mm) and +(down:0mm) .. (OutputVat1.south);		
			
			\node[squarednode] at (7.5,-1.5) [draw=none,fill=none](Pts) {${\color{blue}\mathrm{h}^{\langle k \rangle}}$} ;
			\draw[->,draw=blue,line width=0.25mm] (OutputL1at1.east)  |- ++(0,0) -- (Pts.west);	
			
			\node[squarednode] at (8.2,-3.5) [draw=none,fill=none](InputM) {} ;
			\draw[->,draw=blue,line width=0.25mm] (Pts.south)  ++(0,0) |- (InputM.west);	
			
		\end{tikzpicture}
	}
	\caption{RNN cell computation at time level $k$. $\mathrm{S}$ denotes the sigmoid function.}
	\label{fig:RNN cell iteration}
\end{figure}
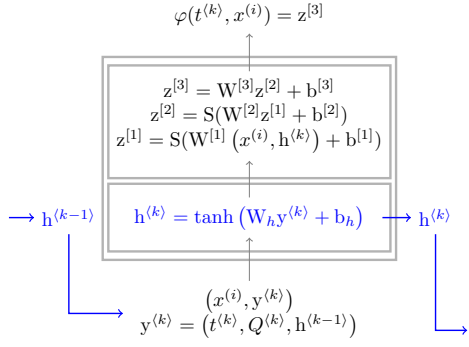

To approximate $\varphi$, 
we consider a RNN. A hidden state, $\mathrm{h}$, carries information about the path history. Thus, the outputs of this architecture depend on the path history and guarantee their progressive measurability.
Our RNN has a hidden layer followed by three dense layers. The cell of the RNN that we iterate over time is depicted in Fig. \ref{fig:RNN cell iteration}; the blue arrows highlight the connection along the temporal component. 

The NN architecture depends on a hyper-parameter $\tau$ than includes all weights $\mathrm{W}$ and biases $\mathrm{b}$ defining the NN. The hyper-parameter $\tau$ is optimized
to decrease  $\mathcal{L}$
by a gradient-descent
method at each training step. 
Let $\varphi_{\Theta}^{\langle k \rangle (i) }$ denote the NN approximation of $\varphi(t_k,x_i)$, where $t_k$ and $x_i$ belong to the time and space grids, respectively. The loss function encodes all the required constraints on the potential. Let 
\begin{align*}
	& \mathcal{L}_{\mathcal{V}}(\varphi_{\tau}) =
	\\
	&  h_x h_t \sum_{k=1}^{M_t} \sum_{i=1}^{M_x} \left\{ F\left(x^{(i)},-\partial_t \varphi_{\tau}^{\langle k \rangle (i) } , \partial_x \varphi_{\tau}^{\langle k \rangle (i) } \right)  \right.
	\\
	& \quad \left.  - u'_T(x^{(i)})\partial_t \varphi_{\tau}^{\langle k \rangle (i) } \right\} ,
	\\
	& \mathcal{L}_{0}(\varphi_{\tau}) = \sum_{k=1}^{M_t} \sum_{i=1}^{M_x} \max\left\{ - \partial_x \varphi_{\tau}^{\langle k \rangle (i) } , 0\right\},
	\\
	& 
	\mathcal{L}_{\mathcal{B}}(\varphi_{\tau}) = \sum_{k=1}^{M_t} \left( h_x \sum_{i=1}^{M_x} \partial_t \varphi_{\tau}^{\langle k \rangle (i) } + Q^{\langle k \rangle} \right)^2 ,
	\\
	& \mathcal{L}_{M_0}(\varphi_{\tau}) = \sum_{i=1}^{M_x} \left(\varphi_{\tau}^{\langle 0 \rangle (i) } - M_0(x^{(i)}) \right)^2 ,
	\\
	& \mathcal{L}_{\mathcal{P}}(\varphi_{\tau}) = \sum_{k=1}^{M_t} \left( 1-h_x\sum_{i=1}^{M_x} \partial_x  \varphi_{\tau}^{\langle k \rangle (i) } \right)^2.
\end{align*}
$ \mathcal{L}_{\mathcal{V}}$ corresponds to the discretization of the variational functional. Deviations from the balance constraint are penalized by $\mathcal{L}_{\mathcal{B}}$, and the initial condition is enforced in $\mathcal{L}_{M_0}$. Moreover, $\mathcal{L}_{0}$ and $\mathcal{L}_{\mathcal{P}}$ guarantee that $\varphi_x(t,\cdot)$ is a density function. Therefore, we consider the loss function
\[
\mathcal{L} =  \mathcal{L}_{\mathcal{V}}+\mathcal{L}_{0}+\mathcal{L}_{\mathcal{B}} +\mathcal{L}_{M_0}+\mathcal{L}_{\mathcal{P}}
\]
that enforces the constraints by penalization.

\subsection{Deterministic price}

At each training step, we evaluate the supply at points on the time grid. Then, for each time level $k$, we compute $\mathrm{h}^{\langle k \rangle}$, which depends on $(t^{\langle i\rangle},Q^{\langle i\rangle})$ for $i\leq k$. Next, we use $\mathrm{h}^{\langle k \rangle}$ for all points $x^{(i)}$ in the state grid. Fig. \ref{fig:RNN cell iteration} illustrates the computation at time level $k$. After iterating over time steps, we obtain the NN approximation of the potential, $\varphi_{\tau}$, at grid points, and we compute $\mathcal{L}$. This process is summarized in Table \ref{tab:Deterministic Algorithm}. 

\begin{table}
	\centering
	\caption{Deterministic training step $j$.}
	\begin{tabular}{|r|l|} 
		\hline
		& Input: $M_0$, $Q$, $\tau^{j-1}$.\\
		\hline 
		{\tiny 1} & for $k=1 ,\ldots , M_t$ do \\
		{\tiny 2} & $\quad \mbox{compute } \mathrm{h}^{\langle k \rangle}$\\
		{\tiny 3} & $\quad \mbox{for } i=1,\ldots , M_x$ do\\
		{\tiny 4} & $\quad \quad \mbox{compute } \varphi_{\tau^{j-1}}(t^{\langle k \rangle},x^{(i)})$\\
		{\tiny 5} & compute $\mathcal{L}(\varphi_{\tau^{j-1}})$\\
		{\tiny 6} & compute $\tau^{j}$ by gradient descend\\
		\hline
		& Output: $\tau^{j}$\\
		\hline
	\end{tabular}
	\label{tab:Deterministic Algorithm}
\end{table}

\subsection{Random price}

The training step follows the same structure as the deterministic case, except for the supply input that is used;  we compute a new random supply sample
for each training step.
The training algorithm is outlined in Table \ref{tab:Stochastic Algorithm}.

\begin{table}
	\centering
	\caption{Stochastic training step $j$.}
	\begin{tabular}{|r|l|} 
		\hline
		& Input: $M_0$, $dQ$, $\tau^{j-1}$.\\
		\hline	 
		{\tiny 1} & Compute sample $Q_j$ \\
		{\tiny 2} & for $k=1 ,\ldots , M_t$ do \\
		{\tiny 3} & $\quad \mbox{compute } \mathrm{h}^{\langle k \rangle}$\\
		{\tiny 4} & $\quad \mbox{for } i=1,\ldots , M_x$ do\\
		{\tiny 5} & $\quad \quad \mbox{compute } \varphi_{\tau^{j-1}}(t^{\langle k \rangle},x^{(i)})$\\
		{\tiny 6} & compute $\mathcal{L}(\varphi_{\tau^{j-1}})$\\
		{\tiny 7} & compute $\tau^{j}$ by gradient descend\\
		\hline
		& Output: $\tau^{j}$\\
		\hline
	\end{tabular}
	\label{tab:Stochastic Algorithm}
\end{table}


\section{Numerical results}
\label{numres}
Here, we take $L(x,v) = \tfrac{1}{2}v^2$ without terminal cost. We select $T=1$ and $M_t = 17$ equally spaced grid points for the time variable. We discretize the space variable interval $[-1,1]$ with $M_x=31$ equally spaced points. 
We set $Q(0)=-0.5$, and $m_0$ is a symmetric and compactly supported distribution within $[-1,1]$. We discretize the partial derivatives as
\begin{align*}
	& \partial_t  \varphi_{\tau}^{\langle k \rangle (i) } = \tfrac{1}{h_t} \left(\varphi_{\tau}^{\langle k+1 \rangle (i) }-\varphi_{\tau}^{\langle k \rangle (i) }\right),
	\\
	& \partial_x  \varphi_{\tau}^{\langle k \rangle (i) } = \tfrac{1}{h_x} \left(\varphi_{\tau}^{\langle k \rangle (i+1) }-\varphi_{\tau}^{\langle k \rangle (i) }\right).
\end{align*}
To guarantee we obtain an approximation of the potential and price on the boundary of $[0,1]\times [-1,1]$, we extend the grid to $[0,1+h_t]\times [-1,1+h_x]$. The ML method is implemented in TensorFlow with Adam optimizer. The dimension of the hidden state is $32$.

\subsection{Deterministic case}
\begin{figure}[htp]
	\centering
	\includegraphics[width=0.35\textwidth]{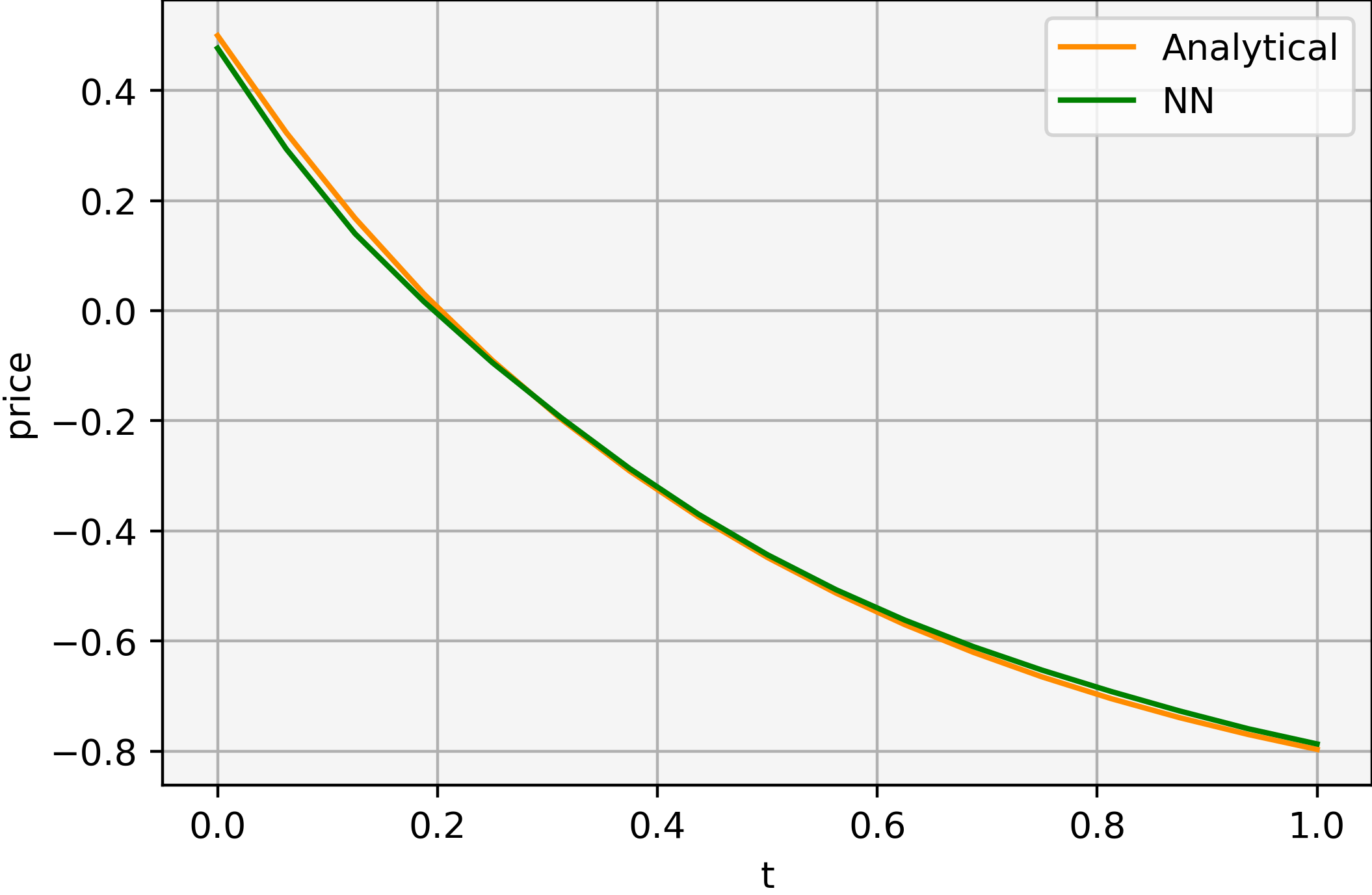}
	\caption{Deterministic price approximation}         
	\label{fig:PricesErrorDeterministic}
\end{figure}

\begin{figure}[htp]
	\centering
	\begin{subfigure}[t]{0.21\textwidth}
		\centering
		\vskip 0pt
		\includegraphics[width=\textwidth]{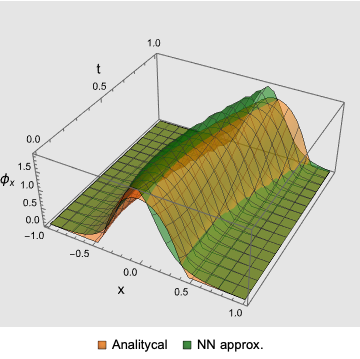}
	\end{subfigure}
	\hfill
	\begin{subfigure}[t]{0.25\textwidth}
		\centering
		\vskip 0pt
		\includegraphics[width=\textwidth]{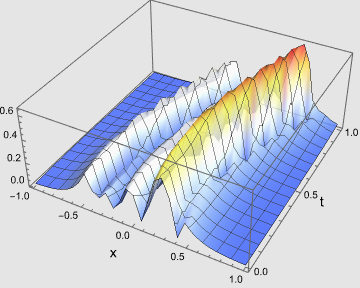}
	\end{subfigure}
	\caption{$m=\varphi_x$ approx. (left) and error (right).}
	\label{fig:Potential_x}
\end{figure}

We assume the supply follows the ODE
\[
dQ(t) = \theta \left( \overline{Q}- Q(t) \right) dt, 
\]
on $[0,T]$, where $\theta = 2$ and $\overline{Q} =1$. In this case, $\varpi = -Q$  \cite{gomes2018mean}. Fig. \ref{fig:PricesErrorDeterministic} shows the price obtained after $18.000$ training steps. The approximation for $m=\varphi_x$
is shown in Figure \ref{fig:Potential_x}. 
For the analytical solution, we have  $\mathcal{L}=\mathcal{L}_{\mathcal{V}}=0.127604$, since $\mathcal{L}_{0}$, $\mathcal{L}_{\mathcal{B}}$, $\mathcal{L}_{M_0}$, and $\mathcal{L}_{\mathcal{P}}$ evaluate to zero. For the NN approximation, we get $\mathcal{L}_{\mathcal{V}}(\varphi_{\tau}) = 0.134984$. These results illustrate an excellent fit for the price function and a good representation of $m$, including its positivity. 

%

\subsection{Stochastic case}

\begin{figure}[htp]
	\centering
	\includegraphics[width=0.35\textwidth]{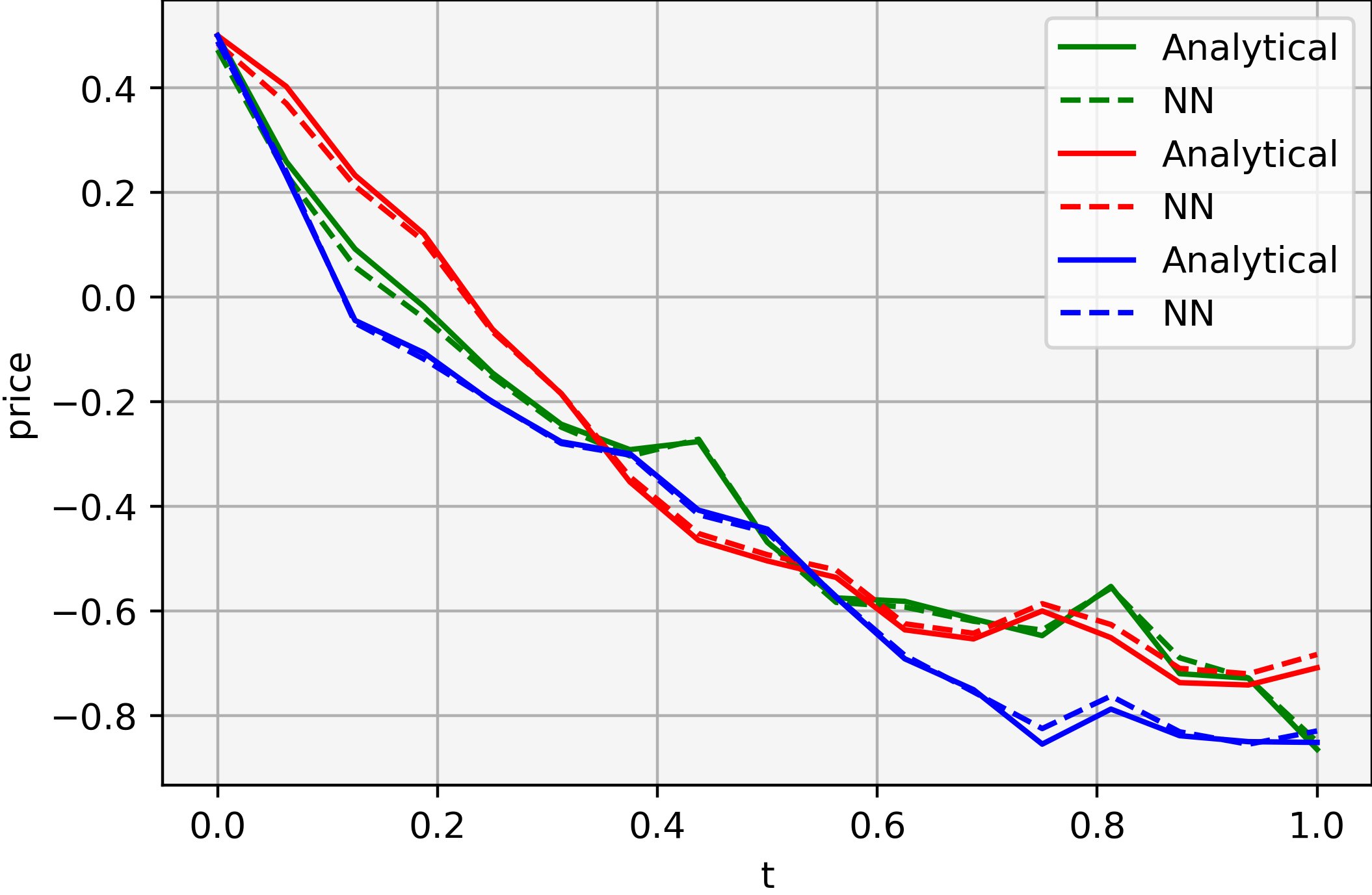}
	\caption{Stochastic price approx. (3 samples)}         
	\label{fig:PricesErrorStochastic}
\end{figure}

\begin{figure}[htp]
	\centering
	\includegraphics[width=0.35\textwidth]{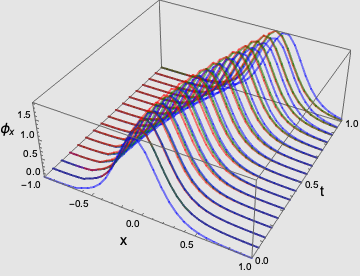}
	\caption{$m=\varphi_x$ approximation (3 samples)}
	\label{fig:PhixStochastic}
\end{figure}

We consider a mean-reverting supply
\[
dQ(t) = \theta \left( \overline{Q}- Q(t) \right) dt+ \sigma dW(t), 
\]
on $[0,T]$, where $\sigma = 0.2$. According to \cite{GoGuRi2021}, the SDE for the price is $d\varpi = -dQ$. To obtain convergence of our algorithm, we increase the training steps to $36000$. We assess the price approximation using a test set of $1000$ samples of the supply, obtaining an $L^\infty$ mean error of $0.036619$. For three samples of the supply, Fig. \ref{fig:PricesErrorStochastic} shows the approximated price, and  Figure \ref{fig:PhixStochastic} shows the approximated $m$. We get
\begin{align*}
	& \mathcal{L}_{\mathcal{V}}(\varphi_{\tau}) = 0.141007, \; \mathcal{L}_{0}(\varphi_{\tau}) = 0, 
	\\
	& \mathcal{L}_{\mathcal{B}}(\varphi_{\tau}) = 0.002728, \; \mathcal{L}_{M_0}(\varphi_{\tau}) = 0.007506, 
	\\
	&\mathcal{L}_{\mathcal{P}}(\varphi_{\tau}) = 0.01905. 
\end{align*}

\section{Conclusions and future work}


We developed a variational approach to price formation with common noise. 
This formulation, combined with machine learning techniques, provides 
an approach to solving infinite-dimensional MFGs. 
Future work should identify better network architectures and 
convergence results. A posteriori estimates would also be 
extremely interesting as they would provide a way to ensure the convergence. 


\section{Acknowledgments}

The authors were supported by King Abdullah University of Science and Technology (KAUST) baseline funds and KAUST OSR-CRG2021-4674.


\bibliographystyle{abbrv}
\bibliography{mfg.bib}

\def\cprime{$'$}
\begin{thebibliography}{10}

\bibitem{CDY}
Y.~Achdou, F.~Camilli, and I.~Capuzzo-Dolcetta.
\newblock Mean field games: numerical methods for the planning problem.
\newblock {\em SIAM J. Control Optim.}, 50(1):77--109, 2012.

\bibitem{DY}
Y.~Achdou and I.~Capuzzo-Dolcetta.
\newblock Mean field games: numerical methods.
\newblock {\em SIAM J. Numer. Anal.}, 48(3):1136--1162, 2010.

\bibitem{MR4214773}
Y.~Achdou, P.~Cardaliaguet, F.~Delarue, A.~Porretta, and F.~Santambrogio.
\newblock {\em Mean field games}, volume 2281 of {\em Lecture Notes in
  Mathematics}.
\newblock Springer, Cham; Centro Internazionale Matematico Estivo (C.I.M.E.),
  Florence, [2020] \copyright 2020.
\newblock Edited by Pierre Cardaliaguet and Alessio Porretta, Fondazione
  CIME/CIME Foundation Subseries.

\bibitem{aid2020equilibrium}
R.~Aid, A.~Cosso, and H.~Pham.
\newblock Equilibrium price in intraday electricity markets, 2020.

\bibitem{AidDumitrescuTankov2021}
R.~A\"{i}d, R.~Dumitrescu, and P.~Tankov.
\newblock The entry and exit game in the electricity markets: A mean-field game
  approach.
\newblock {\em Journal of Dynamics \& Games}, 8(4):331--358, 2021.

\bibitem{ATM19}
C.~Alasseur, I.~Ben~Taher, and A.~Matoussi.
\newblock An extended mean field game for storage in smart grids.
\newblock {\em Journal of Optimization Theory and Applications},
  184(2):644--670, 2020.

\bibitem{alasseur2021mfg}
C.~Alasseur, L.~Campi, R.~Dumitrescu, and J.~Zeng.
\newblock Mfg model with a long-lived penalty at random jump times: application
  to demand side management for electricity contracts, 2021.

\bibitem{SummerCamp2019}
A.~Alharbi, T.~Bakaryan, R.~Cabral, S.~Campi, N.~Christoffersen, P.~Colusso,
  O.~Costa, S.~Duisembay, R.~Ferreira, D.~Gomes, S.~Guo, J.~Gutierrez-Pineda,
  P.~Havor, M.~Mascherpa, S.~Portaro, R.~Ribeiro, F.~Rodriguez, J.~Ruiz,
  F.~Saleh, C.~Strange, T.~Tada, X.~Yang, and Z.~Wr\'{o}blewska.
\newblock A price model with finitely many agents.
\newblock {\em Bulletin of the Portuguese Mathematical Society}, 2019.

\bibitem{YTDJduality2021}
Y.~Ashrafyan, T.~Bakaryan, D.~Gomes, and J.~Gutierrez.
\newblock A duality approach to a price formation mfg model.
\newblock 2021.

\bibitem{YTDJpotential2022}
Y.~Ashrafyan, T.~Bakaryan, D.~Gomes, and J.~Gutierrez.
\newblock A variational approach for price formation models in one dimension,
  2022.

\bibitem{DRT2021Potential}
T.~{Bakaryan}, R.~{Ferreira}, and D.~{Gomes}.
\newblock {A Potential Approach for Planning Mean-Field Games in One Dimension
  }.
\newblock {\em Submitted to Communications on Pure and Applied Analysis}, 2021.

\bibitem{BS02}
T.~{Basar} and R.~{Srikant}.
\newblock Revenue-maximizing pricing and capacity expansion in a many-users
  regime.
\newblock In {\em Proceedings.Twenty-First Annual Joint Conference of the IEEE
  Computer and Communications Societies}, volume~1, pages 294--301 vol.1, 2002.

\bibitem{Csato2011ThePE}
G.~Csat'o, B.~Dacorogna, and O.~Kneuss.
\newblock {\em The Pullback Equation for Differential Forms}.
\newblock Progress in Nonlinear Differential Equations and their Applications.
  Birkh\"{a}user/Springer, New York, 2012.

\bibitem{TBD20}
B.~Djehiche, J.~Barreiro-Gomez, and H.~Tembine.
\newblock {Price Dynamics for Electricity in Smart Grid Via Mean-Field-Type
  Games}.
\newblock {\em Dynamic Games and Applications}, 10(4):798--818, December 2020.

\bibitem{FTT20}
O.~F{\'e}ron, P.~Tankov, and L.~Tinsi.
\newblock {Price Formation and Optimal Trading in Intraday Electricity Markets
  with a Major Player}.
\newblock {\em Risks}, 8(4):1--1, December 2020.

\bibitem{feron2021price}
O.~F\'{e}ron, P.~Tankov, and L.~Tinsi.
\newblock Price formation and optimal trading in intraday electricity markets,
  2021.

\bibitem{FT20}
M.~Fujii and A.~Takahashi.
\newblock {A Mean Field Game Approach to Equilibrium Pricing with Market
  Clearing Condition}.
\newblock Papers 2003.03035, arXiv.org, Mar. 2020.

\bibitem{fujii2021equilibrium}
M.~Fujii and A.~Takahashi.
\newblock Equilibrium price formation with a major player and its mean field
  limit, 2021.

\bibitem{GoGuRi2021B}
D.~{Gomes}, J.~{Gutierrez}, and R.~{Ribeiro}.
\newblock {A random-supply Mean Field Game price model}.
\newblock {\em arXiv e-prints}, page arXiv:2109.01478, Sept. 2021.

\bibitem{GoGuRi2021}
D.~Gomes, J.~Gutierrez, and R.~Ribeiro.
\newblock A mean field game price model with noise.
\newblock {\em Math. Eng.}, 3(4):Paper No. 028, 14, 2021.

\bibitem{gomes2018mean}
D.~Gomes and J.~Sa\'{u}de.
\newblock A {M}ean-{F}ield {G}ame {A}pproach to {P}rice {F}ormation.
\newblock {\em Dyn. Games Appl.}, 11(1):29--53, 2021.

\bibitem{Caines1}
M.~Huang, R.~P. Malham{\'e}, and P.~E. Caines.
\newblock Large population stochastic dynamic games: closed-loop
  {M}c{K}ean-{V}lasov systems and the {N}ash certainty equivalence principle.
\newblock {\em Commun. Inf. Syst.}, 6(3):221--251, 2006.

\bibitem{ll1}
J.-M. Lasry and P.-L. Lions.
\newblock Jeux \`a champ moyen. {I}. {L}e cas stationnaire.
\newblock {\em C. R. Math. Acad. Sci. Paris}, 343(9):619--625, 2006.

\bibitem{BS10}
H.~Shen and T.~Basar.
\newblock Pricing under information asymmetry for a large population of users.
\newblock {\em Telecommun. Syst.}, 47(1-2):123--136, 2011.

\bibitem{JSF20}
A.~Shrivats, D.~Firoozi, and S.~Jaimungal.
\newblock {A Mean-Field Game Approach to Equilibrium Pricing, Optimal
  Generation, and Trading in Solar Renewable Energy Certificate Markets}.
\newblock Papers 2003.04938, arXiv.org, Mar. 2020.

\end{thebibliography}

\end{document}